\documentclass{article}

\usepackage{amsmath}
\usepackage{amsfonts}
\usepackage{amssymb}
\usepackage{mathrsfs}
\usepackage{graphicx}

\newcommand   \Integers {\mathbb Z}
\newcommand   \PositiveIntegers {{\mathbb N}}
\newcommand   \primeInt {\mathtt{q}}
\newcommand   \fieldChar {\mathtt{\,p\,}}

\renewcommand   \mod  {~\textsf{mod}~ }
\renewcommand   \gcd  {~\textsf{gcd}~ }

\newcommand  \basel[2]{#1_{_{#2}}}

\newcommand  \lspace  {\hspace*{-1.5mm}}

\newcommand  \tab  {\hspace*{0.5cm}}
\newcommand  \ltab  {\hspace*{-0.5cm}}

\newcommand  \shiftright  {\hspace*{2.0cm}}
\newcommand  \shiftleft  {\hspace*{-2.0cm}}

\newtheorem{proposition}{Proposition}

\newcommand \proof {\rm{\textbf{Proof.}}~~}

\newcommand {\qed} {\hfill{} $\square$}

\newcommand  \bglb {\big (}
\newcommand  \bgrb {\big )}

\newcommand  \bgls {\big [}
\newcommand  \bgrs {\big ]}

\begin{document}

\title{Probability Analysis and Comparison of Well-Known Integer Factorization Algorithms}

\author{\bf{Duggirala Meher Krishna}\\
{\small{Gayatri Vidya Parishad College of Engineering (Autonomous)}} \\
{\small{Madhurawada, VISAKHAPATNAM -- 530 048, Andhra Pradesh, India}} \\
 {\small{E-mail ~: \tab duggiralameherkrishna@gmail.com}}\\
 \\
 and \\
 \\
\bf{Duggirala Ravi}\\
{\small{Gayatri Vidya Parishad College of Engineering (Autonomous)}} \\
{\small{Madhurawada, VISAKHAPATNAM -- 530 048, Andhra Pradesh, India}} \\
\shiftleft \tab {\small{E-mail ~: \tab ravi@gvpce.ac.in; \tab duggirala.ravi@yahoo.com};} \\
{\small{\shiftright  duggirala.ravi@rediffmail.com; \tab drdravi2000@yahoo.com}} 
}

\date{}

\maketitle

\begin{abstract}
Two prominent methods for integer factorization are those based on general integer sieve and elliptic curve. The general integer sieve method can be specialized to quadratic integer sieve method. In this paper, a probability analysis for the success of these methods is described, under some reasonable conditions. The estimates presented are specialized for the elliptic curve factorization. These methods are compared through heuristic estimates. It is shown that the elliptic curve method is a probabilistic polynomial time algorithm under the assumption of uniform probability distribution for the arising group orders and clearly more likely to succeed, faster asymptotically. 
\\

\\


\end{abstract}

\section{Introduction}
In this paper, the success probabilities for two prominent methods, {\em viz}, general integer sieve method and elliptic curve method, are presented. The estimates are specialized for the elliptic curve factorization algorithm. The random variables studied are (1) the number generated by exponentiating a chosen fixed base random number to various random integer exponents, for general integer sieve method, and (2) the group orders of the elliptic curve groups, with restriction to $\mod \fieldChar$, for each (as yet unknown) prime factor $\fieldChar$ of the integer modulus to be factored. The common assumptions taken in our estimates are that the probabilistic events arising from the consideration of various different smaller prime numbers being factors of any particular realization (sample) of the random variable are mutually independent. With the assumption of independence of events corresponding to divisibility by different smaller prime numbers, the probabilities of success are shown to be fairly optimistic. The general integer sieve needs the random base point to be a group generator (primitive in this sense), which may be difficult to ensure. The merits of elliptic curve method are highlighted, with a caution concerning the widths of the intervals of the possible group orders. Nevertheless, the estimated probabilities of success do not depend too heavily on this fact, as they are applicable to random samples form any arbitrary interval of considerable width, for asymptotic analysis.

  
\section{\label{Sec-Formulation}Estimation of Success Probabilities}
Let $\Integers$ be the ring of integers, and $\PositiveIntegers$ be the set of positive integers.
Let $N$ be a very large positive integer to be factored, and let $\basel{\Integers}{N}$ be the 
ring of integers with arithmetic operations taken $\mod N$.

 Let $L_{\min}, L_{\max} \in \Integers$ be such that $L_{\min} < L_{\max}$ and $L_{\max}-L_{\min}$ is very large. The consecutive prime numbers are listed in the ascending order as follows: $2 = \basel{\primeInt}{1}, \, 3 = \basel{\primeInt}{2}, \, 5 = \basel{\primeInt}{3}, ....$, so that $\basel{\primeInt}{i}$ is the $i$-th prime number, for $i \in \PositiveIntegers$.  Let $k$ be a small positive integer, but still large enough that the asymptotic estimates hold good,  and let $n$ be the largest positive integer, such that $\basel{\primeInt}{n} < \max \{ |L_{\min}|, \, |L_{\max}| \}$. Let $X$ be a random variable taking integer values in the interval ${\mathcal I} =  \bgls L_{\min}\, ,~ L_{\max}\bgrs$, with uniform probability distribution. 

\begin{proposition}
\label{prop-01}
In the notation just discussed, the probability $\basel{\pi}{X}(z)$ of the event that
a sample of the random variable $X$ is divisible by a positive integer $z \geq 2$
is approximately $\frac{1}{z}$, and more precisely the following bounds hold good:
       \begin{equation}
           \frac{1}{z} - \frac{1}{L_{\max}-L_{\min}} ~~ \leq ~~ \basel{\pi}{X}(z)
               ~~ \leq ~~       \frac{1}{z} + \frac{1}{L_{\max}-L_{\min}}  
                          \label{Bounds-on-pi-X-z}
       \end{equation}
\end{proposition}
\proof For every positive integer $ z \geq 2$, the number of integer multiples of $z$ in ${\mathcal I}$ are between
$\bglb \frac{L_{\max}-L_{\min}}{z}-1 \bgrb$ and $\bglb \frac{L_{\max}-L_{\min}}{z}+1 \bgrb$.
Thus, the probability that a random sample of $X$ is divisible by $z$ is between
$\frac{1}{z} - \frac{1}{L_{\max}-L_{\min}}$ and $\frac{1}{z} +  \frac{1}{L_{\max}-L_{\min}}$,
which justifies the assumptions, with appropriate choices of $z$.     \qed\\

The conjunct consideration concerning the divergence of $ \sum_{i} \frac{1}{\basel{\primeInt}{i}}$
and the convergence of $ \sum_{i} \frac{1}{\basel{\primeInt^{2}}{i}}$ necessitates taking product spaces.
Moreover, the estimates are presented only for elliptic curve factorization algorithm.

\subsection{\label{Sec-Elliptic-Curve}Success of Elliptic Curve Factorization}
Let $r = \left \lceil \frac{\log (N)}{\log (\basel{\primeInt}{k})} \right \rceil$, 
where the choice of $k$, the number of smaller prime factors to be used,
is assumed to be considerably larger than $2$, such as about $1000$. Actually,
$k$ can run into tens of thousands, for practical purposes, and constrained
by the condition that $\basel{\primeInt^{r}}{k} \geq N$.
If $\basel{\primeInt}{k}$ is too small, then $r$ can be so large
that the estimated failure probabilities may become irrelevant.
Let $\basel{\mathcal C}{l}\bglb \basel{\Integers}{N} \bgrb$ be elliptic curves,
defined over $\basel{\Integers}{N}$, for $1 \leq l \leq r$.
Let $\fieldChar$ be a large but unknown prime integer factor $N$, 
such that $\fieldChar \leq \sqrt{N}$, and 
$\basel{\mathcal C}{l}\bglb\basel{\Integers}{\fieldChar}\bgrb$
be the corresponding elliptic curves restricted to $\basel{\Integers}{\fieldChar}$,
for $1 \leq l \leq r$.  The group order of 
$\basel{\mathcal C}{l}\bglb\basel{\Integers}{\fieldChar}\bgrb$
is $\fieldChar+1-\basel{a}{l}$, where $-2\sqrt{\fieldChar} \leq \basel{a}{l} \leq 2\sqrt{\fieldChar}$,
by Hasse-Weil bounds for the elliptic curve group orders.
The probability distribution of $\fieldChar+1-t$ of the group order
of ${\mathcal C}\bglb\basel{\Integers}{\fieldChar}\bgrb$,
as obtained by taking $\mod \fieldChar$ restriction of a randomly
generated elliptic curve 
${\mathcal C}\bglb \basel{\Integers}{N} \bgrb$
is assumed to be uniform over the interval
${\mathcal I} = [(\sqrt{\fieldChar}-1)^2\, ,~  (\sqrt{\fieldChar}+1)^2]$.

\begin{proposition}
\label{prop-02}
Let $\basel{\mathcal C }{l}\bglb \basel{\Integers}{N} \bgrb$, 
for $1 \leq l \leq r+2$, be any $(r+2)$ independent samples of the
elliptic curves, and $\fieldChar$ be a fixed (though unknown yet)
prime factor of $N$, such that $\fieldChar \leq \sqrt{N}$.
Let ${\mathcal E}_{k+1}$ be the random event that each of the $(r+2)$
group orders  $\fieldChar + 1 - \basel{a}{l}$ of the elliptic curves
 $\basel{\mathcal C }{l}\bglb \basel{\Integers}{\fieldChar} \bgrb$, 
for $1 \leq l \leq r+2$, is divisible by a prime factor
at least as large as $\basel{\primeInt}{k+1}$, where the prime
number $\fieldChar$ is assumed to be such that $\fieldChar \, \vert \, N$
and $\fieldChar \geq \basel{\primeInt}{k+1}$. Then,
$Pr\bglb {\mathcal E}_{k+1}\bgrb \leq \frac{(r+2)(r+1)+8}{2 \times 4 \times (\basel{\primeInt}{k+1}-1)}$
$~+~$
${\mathcal O}\bglb \frac{(r+2)(r+1)}{8}\times \frac{ \log (~ \log (\fieldChar) ~ )}{\sqrt{\fieldChar}}\bgrb$.
 Further, if the approximation $\basel{q}{i} \approx i \log (i)$, for sufficiently
 large positive integer $i$,  is permitted, then
 $Pr\bglb {\mathcal E}_{k+1}\bgrb \leq \frac{(r+2)(r+1)+8}{2 \times 4 \times k \times (\log (k+1))^{2}}$
 $~+~$
${\mathcal O}\bglb \frac{(r+2)(r+1)}{8}\times \frac{ \log (~ \log (\fieldChar) ~ )}{\sqrt{\fieldChar}}\bgrb$.
\end{proposition}
\proof Before proceeding with the proof, a justification for the
validity of the approximation  in the last part is as follows:
by the prime number theorem,
$i \approx \frac{\basel{\primeInt}{i}}{\log(\basel{\primeInt}{i})} < \frac{\basel{\primeInt}{i}}{\log(i)}$,
and $ \basel{\primeInt}{i} $ is likely to be larger than $i \log (i)$.
It may also be noticed that 
$\frac{(r+2)(r+1)+8}{8k (\log (k+1))^{2}} \approx \frac{(r+2)(r+1)+8}{8\basel{\primeInt}{k} (\log (k+1))}$.

The random event ${\mathcal E}_{k+1}$ in the statement is broken up into the following two parts:
${\mathcal E}_{k+1} \subseteq E_{k+1, \, 1} \cup E_{k+1, \, 2}$, where

\begin{enumerate}
\item $E_{k+1, \, 1}$ is the event that there are distinct prime numbers 
$\basel{\primeInt}{\basel{i}{l}} \geq \basel{\primeInt}{k+1}$,
for $1 \leq l \leq r+2$, such that $\basel{\primeInt}{\basel{i}{l}} \, \mid \, (\fieldChar+1-\basel{a}{l})$ and
$\basel{\primeInt}{\basel{i}{l}} \, \nmid \, (\fieldChar+1-\basel{a}{l'})$, for $l' \neq l$ and
$1 \leq l,\, l \leq r+2$,   and 

\item  $E_{k+1, \, 2}$ is the event that there is a prime number
$\basel{\primeInt}{i} \geq \basel{\primeInt}{k+1}$, such that
$\basel{\primeInt}{i} \, \mid \, (\fieldChar+1-\basel{a}{l})$ and
$\basel{\primeInt}{i} \, \mid \, (\fieldChar+1-\basel{a}{l'})$,
for two indexes $l$ and $l'$, $l' \neq l$, where $1 \leq l,\, l \leq r+2$.

\end{enumerate}

The two events listed above are not mutually exclusive, but an upper found for
the sum of their probabilities is found, as an estimate for the upper bound
of the event in the statement.

\paragraph{Part (1). ~~}For the event $E_{k+1,\, 1}$,
 it is observed that,  from the simultaneous congruence relations
 $\fieldChar +1 \equiv \basel{a}{l} \mod \basel{\primeInt}{\basel{i}{l}}$,
 for $1 \leq l \leq r+2$, the fixed number $\fieldChar+1$ can be recovered
 by the Chinese remainder theorem. The mapping 
 $\basel{a}{l} \mapsto \basel{a}{l} \mod \basel{\primeInt}{\basel{i}{l}}$,
 for $1 \leq l \leq r+2$,  induces the homomorphism 
 $(\basel{a}{1}, \, \cdots, \, \basel{a}{r+2}) \mapsto 
 (\basel{a}{1} \mod \basel{\primeInt}{\basel{i}{1}}, \, \cdots, \, \basel{a}{r+2} \mod \basel{\primeInt}{\basel{i}{r+2}})$,
 that preserves the algebraic structure. In the proof, it is assumed that the probability distributions remain
 uniform under the mapping $\basel{a}{l} \mapsto \basel{a}{l} \mod \basel{\primeInt}{\basel{i}{l}}$,
 for $1 \leq l \leq r+2$,  with restriction on the domain of possible values of 
 $ (\basel{a}{1} \mod \basel{\primeInt}{\basel{i}{1}}, \, \cdots, \, \basel{a}{r+2} \mod \basel{\primeInt}{\basel{i}{r+2}})$.

 By the mutual independence of $\basel{a}{l}$, for $1 \leq i \leq r+2$, 
 there are at least $4^{r+2} \prod_{l = 1}^{r+2} \sqrt{\basel{\primeInt}{\basel{i}{l}}}$ many possibilities,
 in all, for the set of possible realizations
 $ (\basel{a}{1} \mod \basel{\primeInt}{\basel{i}{1}}, \, \cdots, \, \basel{a}{r+2} \mod \basel{\primeInt}{\basel{i}{r+2}})$, after taking into account the restriction that $|a_{l}| \leq 2\sqrt{\fieldChar}$.
 The fixed number  $\fieldChar+1$ must belong to the set of positive integers
 that can be reconstructed by any realization of
  $ (\basel{a}{1} \mod \basel{\primeInt}{\basel{i}{1}}, \, \cdots, \, \basel{a}{r+2} \mod \basel{\primeInt}{\basel{i}{r+2}})$, with $\fieldChar$ constrained to be a prime number.
  Now, the number of possibilities for the realizations
  for  $ (\basel{a}{1} \mod \basel{\primeInt}{\basel{i}{1}}, \, \cdots, \, \basel{a}{r+2} \mod \basel{\primeInt}{\basel{i}{r+2}})$, that could result in the reconstruction of $\fieldChar+1$,
  with $\fieldChar$ restricted to be a prime number at most $\sqrt{N}$
    (or of bit size at most $\frac{\log_{2}(N)}{2}$),
  is smaller than   $\prod_{l = 1}^{r} \sqrt{\basel{\primeInt}{\basel{i}{l}}}$,
  because $\bglb\sqrt{\basel{\primeInt}{k}}\bgrb^{r} \geq \sqrt{N} > \frac{\fieldChar+1}{2}$.
  Thus, $Pr\bglb E_{k+1,\, 1} \bgrb$
  $ \leq $
  $\frac{1}{\sqrt{\basel{\primeInt}{\basel{i}{r+1}}\basel{\primeInt}{\basel{i}{r+2}}}}$
   $\leq \frac{1}{\basel{\primeInt}{k+1}}$.
   A justification for this approach is given in a separate
   paragraph following the proof of the second part.
 
 \paragraph{Part (2). ~~} For the event $E_{k+1,\, 2}$,
 a slightly weaker proof is given in this paragraph,
 and a more accurate proof is given the correction part below.
 The event that a prime number
 $\basel{\primeInt}{i} \geq \basel{\primeInt}{k+1}$,
 such that $\basel{\primeInt}{i}$ divides the group orders of both
 $\basel{\mathcal C}{l}\bglb \basel{\Integers}{N}\bgrb$ and
  $\basel{\mathcal C}{l'}\bglb \basel{\Integers}{N}\bgrb$, for some $l$ and $l'$,
 $l \neq l'$ and $1 \leq l, \, l' \leq r+2$, occurs with probability
 $\frac{(r+2)(r+1)}{2 \basel{\primeInt^{2}}{i}}$, for any $i$, where $i \geq k+1$.
 This probability also accounts for the possibility that
 $\basel{\primeInt}{i} \, \vert \, \fieldChar+1-\basel{a}{l}$
 and  $\basel{\primeInt}{i} \, \vert \, \fieldChar+1-\basel{a}{l'}$,
 in case $\basel{a}{l} = \basel{a}{l'}$, but $l \neq l'$,
 where $1 \leq l, \, l' \leq r+2$, for some prime number
 $\fieldChar \, \mid \, N$ and $\fieldChar \geq \basel{\primeInt}{k+1}$.
 However, there are at least four possibilities that $\basel{\primeInt}{i}$
 divides either component of the pairs
   $(\fieldChar+1-\basel{a}{l} \, , \, \fieldChar+1-\basel{a}{l'})$,
     $(\fieldChar'+1-\basel{a'}{l} \, , \, \fieldChar+1-\basel{a}{l'})$,
          $(\fieldChar+1-\basel{a}{l} \, , \, \fieldChar'+1-\basel{a'}{l'})$
  and      $(\fieldChar'+1-\basel{a'}{l} \, , \, \fieldChar'+1-\basel{a'}{l'})$,
  for  two distinct prime factors $\fieldChar$ and $\fieldChar'$ of the composite
 number $N$, of which only one possibility is taken into account, for a fixed $\fieldChar$.
  Thus, a multiplier by at most the fraction $\frac{1}{4}$ must be applied. 
  Now, $\sum_{i \geq k+1} \frac{1}{\basel{\primeInt^{2}}{i}} < $
 $\sum_{i \geq k+1} \bgls \frac{1}{\basel{\primeInt}{i}-1}-\frac{1}{\basel{\primeInt}{i}} \bgrs$.
 $<  \frac{1}{\basel{\primeInt}{k+1}-1}$. The result follows by adding it to probability bound in 
 the first part. 
 
   If the approximation $\basel{q}{i} \approx i \log (i)$ is permitted, the probability bound
   in the second part is as follow: 
   $\sum_{i \geq k+1} \frac{1}{\basel{\primeInt^{2}}{i}} \approx $
   $\sum_{i \geq k+1} \frac{1}{i^{2} (\log (i))^{2}} < $ 
   $\frac{1}{(\log (k+1))^{2}}\sum_{i \geq k+1} \frac{1}{i^{2}} < $ 
   $\frac{1}{(\log (k+1))^{2}}\sum_{i \geq k+1} \bgls \frac{1}{i-1} - \frac{1}{i} \bgrs $ 
   $ < \frac{1}{k(\log (k+1))^{2}}$.   \qed
 \\
 
 In the following, a justification for the upper bound for $Pr\bglb E_{k+1, \, 1} \bgrb$
 and a small correction to the upper bound for $Pr\bglb E_{k+1,\, 2} \bgrb$,
 assuming that $N$ is a random integer modulus of a prescribed bit size, are given.
 
 \paragraph{Justification for Upper Bound for $Pr\bglb E_{k+1,\, 1} \bgrb$.~~}
  Conditional and joint probabilities over the possible random modulus integer $N$,
  of bit size equal to a prescribed parameter
  $(\lceil \basel{\log}{2}(N)\rceil)$, for independent realizations 
  of the tuples $(\basel{a}{1}, \, \ldots, \, \basel{a}{r+2})$,
  with appropriate restrictions on the domains of possible values,
  are taken into consideration. 
 Let the sequences $(\basel{i}{1}, \, \ldots, \basel{i}{r+2})$,  for 
  $\basel{i}{l} \neq \basel{i}{l'}$ and $k+1 \leq \basel{i}{l},\, \basel{i}{l'} \leq n$,
  where $1 \leq l, \, l' \leq r+2$, $l \neq l'$ and 
  $n$ is the largest positive integer such that
  $\basel{q}{n} \leq (N^{\frac{1}{4}}+1)^{2}$, be enumerated in some particular total order,
  denoted by $\prec$.
  Let $X_{(\basel{i}{1}, \, \ldots, \basel{i}{r+2})}$ be the event that the group order of
  $\basel{\mathcal C }{l}\bglb \basel{\Integers}{N} \bgrb$ is divisible by $\basel{\primeInt}{\basel{i}{l}}$,
  for $1 \leq l \leq r+2$, over all possible integer moduli of bit size $(\lceil \basel{\log}{2}(N)\rceil)$,
  excluding the events $X_{(\basel{j}{1}, \, \ldots, \basel{j}{r+2})}$, for 
  $(\basel{j}{1}, \, \ldots, \basel{j}{r+2}) \prec (\basel{i}{1}, \, \ldots, \basel{i}{r+2})$,
  if any. Now
  \begin{small}
  \begin{eqnarray*}
&& \ltab \ltab  Pr\bglb E_{k+1,\, 1} \bgrb ~~ \leq ~~
   \sum_{(\basel{i}{1}, \, \ldots, \basel{i}{r+2})}  \bgls \tab \tab   Pr\bglb X_{(\basel{i}{1}, \, \ldots, \basel{i}{r+2})} \bgrb  ~ \times \\
&& \tab \tab \tab \tab \tab \tab   Pr\bglb \tab  \textrm {the event that~} \fieldChar \textrm{~ is a large prime number} \\
&&  \tab \tab \tab \tab \tab \tab \tab \tab \textrm{~~ of bit size at most~} \frac{\log_{2}(N)}{2} , \textrm{~such that,}\\
&&  \tab \tab \tab \tab \tab \tab  \tab \tab \textrm{~~ for every~} l, ~~ \basel{\primeInt}{\basel{i}{l}} \, \mid \, \fieldChar + 1 - \basel{a}{l} ,  \textrm{~~and} \\
 &&  \tab \tab \tab \tab \tab \tab \tab \textrm{~~ for some~} l', ~~ \basel{\primeInt}{\basel{j}{l'}} \, \nmid \, \fieldChar + 1 - \basel{a}{l'} ,  \textrm{~ whenever~}  \\
 && \tab \tab \tab \tab \tab \tab \tab \tab \tab \tab (\basel{j}{1}, \, \ldots, \basel{j}{r+2})\prec (\basel{i}{1}, \, \ldots, \basel{i}{r+2})\,,  \\
  && \tab \tab \tab \tab \tab \tab \tab \tab \tab   \textrm{~where~} 1 \leq l,\, l' \leq r+2 \tab \bgrb \tab \tab \bgrs\\
&& \tab \tab  \leq ~~ 
   \sum_{(\basel{i}{1}, \, \ldots, \basel{i}{r+2})}  \bgls \tab \tab   Pr\bglb X_{(\basel{i}{1}, \, \ldots, \basel{i}{r+2})} \bgrb  ~ \times \\
&& \tab \tab \tab \tab \tab \tab   Pr\bglb \tab  \textrm {the event that~} \fieldChar \textrm{~ is a large prime number} \\
&&  \tab \tab \tab \tab \tab \tab \tab \tab \textrm{~~ of bit size at most~} \frac{\log_{2}(N)}{2} , \textrm{~such that,}\\
&&  \tab \tab \tab \tab \tab \tab  \tab \tab \textrm{~~ for every~} l, ~~ \basel{\primeInt}{\basel{i}{l}} \, \mid \, \fieldChar + 1 - \basel{a}{l} , \\
 && \tab \tab \tab \tab \tab \tab \tab \tab \tab   \textrm{~where~} 1 \leq l \leq r+2 \tab \bgrb \tab \tab \bgrs\\
&& \tab \tab \leq ~~ \sum_{(\basel{i}{1}, \, \ldots, \basel{i}{r+2})}   Pr\bglb X_{(\basel{i}{1}, \, \ldots, \basel{i}{r+2})} \bgrb  ~  \times ~ \frac{1}{\basel{\primeInt}{k+1}} ~~~~ \leq ~~~~ \frac{1}{\basel{\primeInt}{k+1}}
  \end{eqnarray*}
  \end{small}
  
\paragraph{Small Correction of Upper Bound for $Pr\bglb E_{k+1,\, 2} \bgrb$. ~~} Taking the upper estimate
$ \frac{1}{\basel{\primeInt}{i}}+\frac{1}{4\sqrt{\fieldChar}}$ in place of
$\frac{1}{\basel{\primeInt}{i}}$, for $k+1 \leq i \leq n$, the following is obtained:
\begin{small}
\begin{eqnarray*}
  && \ltab       Pr\bglb E_{k+1, 2} \bgrb ~~ \leq ~~ \sum_{i = k+1}^{n} \left ( \frac{1}{\basel{\primeInt}{i}}+\frac{1}{4\sqrt{\fieldChar}}\right )^{2} 
  ~~ = ~~ \sum_{i = k+1}^{n} \left ( \frac{1}{\basel{\primeInt^{2}}{i}} ~ + ~
  \frac{1}{8 \basel{\primeInt}{i} \sqrt{\fieldChar}}   ~ + ~\frac{1}{16\fieldChar}\right ) 
  \end{eqnarray*}
 \end{small}
 \lspace where $n$ is constrained to be the largest positive integer such that
  $\basel{\primeInt}{n}$ may possibly divide both $\fieldChar+1-a$ and
  $\fieldChar+1-a'$,  for some $-2\sqrt{\fieldChar} \leq a,\, a' \leq 2\sqrt{\fieldChar}$.
  Since $\gcd(\fieldChar+1-a\, ,\, \, \fieldChar+1-a')$ must divide $|a-a'| \leq 4\sqrt{\fieldChar}$,
  it may be assumed that $n \leq \frac{4\sqrt{\fieldChar}}{\log (4\sqrt{\fieldChar})}$,
  when $a \neq a'$. The terms accrued from
  \begin{small}
  \begin{enumerate}
  \item the sum $\frac{1}{\sqrt{\fieldChar}}\sum_{i = k+1}^{n} \frac{1}{\basel{\primeInt}{i}}$,
  which can be replaced with
 $\frac{\log \bglb \log (\basel{\primeInt}{n}) \bgrb}{\sqrt{\fieldChar}} \approx \frac{\log \bglb 2\log (\sqrt{\fieldChar}+1) \bgrb}{\sqrt{\fieldChar}}~$; 
   \item the event that $a = a'$, which is $\frac{1}{4\sqrt{\fieldChar}}$,
  for independent samples $a$ and $a'$, assuming values from the interval
  $[-2\sqrt{\fieldChar}\, ,\, \, 2\sqrt{\fieldChar}]$ ;  and
  \item the sum $\sum_{i = k+1}^{n} \frac{1}{\fieldChar}$,
  which can be replaced with
 $ \frac{(4\sqrt{\fieldChar})}{ \fieldChar\log (4\sqrt{\fieldChar}) }$
 $ = $ 
 $ \frac{4}{\sqrt{\fieldChar} \log (4\sqrt{\fieldChar}) }$ 
 \end{enumerate}
 \end{small}
 are insignificant for large $\fieldChar$.
 In the statement of the proposition, the effect of the correction terms 
 is reflected in the addend 
${\mathcal O}\bglb \frac{(r+2)(r+1)}{8}\times \frac{ \log (~ \log (\fieldChar) ~ )}{\sqrt{\fieldChar}}\bgrb$.
 \\
 
 The methods for justification and correction terms are similar to
 {\em{a priori}} and {\em{a posteriori}} estimation of the probabilities.
 To be more explicit, the probability that a random prime $\fieldChar$
 being a factor of the random modulus $N$, where $N$ satisfies the
 requirements specified by  $X_{(\basel{i}{1}, \, \ldots, \basel{i}{r+2})}$,
 with specified bit size of $\basel{\log}{2}(N)$ of a fixed number, 
 assuming uniform likelihood among all such prime numbers
 that may arise, is estimated and shown to be upper bounded by 
 $\frac{1}{\basel{\primeInt}{k+1}}$. If we were to take $\frac{1}{\fieldChar}$
 for the probability distribution of this event, we would, actually,
 get an even smaller upper bound for  $Pr\bglb E_{k+1,\, 1} \bgrb$.
 This indirect approach is necessitated by the difficulties
 arising out of the need to deal with the principle of
 inclusion-and-exclusion in the estimation of the  probability
 of union of events, from the probabilities of
 independent individual atomic events.
 For instance, if $Pr\bglb {\mathcal E}_{k+1} \bgrb$ is replaced with something like
 $\frac{\sum_{i = k+1}^{n} \frac{1}{\basel{\primeInt}{i}}}{\sum_{i = 1}^{n} \frac{1}{\basel{\primeInt}{i}}}$,
 for some large enough $n$, the resulting failure probability may become totally unrealistic.
 If hyperelliptic curve method can be adapted for factorization,
 the success probability may hopefully become better.

\section{\label{Sec-General-Comparison}Comparison with General Integer Sieve Factorization}
Let $N$ be a large composite positive integer, and $g \in \basel{\Integers^{*}}{N}$,
where $\basel{\Integers^{*}}{N}$ is the group of invertible elements $\mod N$,
with respect to the multiplication $\mod N$. For a randomly chosen $t \in \basel{\Integers}{N}$,
estimates for the probability of the event that every prime factor of $g^{t} \mod N$ is at most
$\basel{\primeInt}{k}$ remain elusive.  The operational theory of general integer sieve method
is described below.

Let $d_{j}$ be the discrete logarithm of $\basel{\primeInt}{j}$, assuming that $\basel{\primeInt}{j}$
belong to the cyclic subgroup generated by $g$, for $1 \leq j \leq k$.
After collecting sufficient number of samples, a system linear equations of the form
$\sum_{j = 1}^{k} \basel{\nu}{i, \, j} d_{j} \equiv \basel{t}{i} \mod \phi(N)$ is formed, for $1 \leq i \leq k$,
where $\phi(N)$ is the Euler function of $N$, which is the group order of $ \basel{\Integers^{*}}{N}$.
Any such relation arise as a result of the factorization 
 $g^{\basel{t}{i}} = \prod_{j = 1}^{k} \basel{\primeInt^{\basel{\nu}{i, j}}}{i}$,
 for some random samples $\basel{t}{i}$, for $1 \leq i \leq k$.
 
From every new relation 
$\sum_{j = 1}^{k} \basel{\nu}{k+l, \, j} d_{j} \equiv \basel{t}{k+l} \mod \phi(N)$,
a vector, consisting of integers $\basel{\tau}{k+l, \, i}$, $1 \leq i \leq k$,
as components, may be hopefully found, such that
$\sum_{i = 1}^{k} \basel{\tau}{k+l, \, i}\basel{\nu}{i, \, j}  \equiv 0 \mod \phi(N)$,
for $l = 1, 2, 3, \ldots$. Some of the relations may be redundant, leading to trivial relations.
In fact, if two linearly independent relations 
$\sum_{j = 1}^{k} \basel{\nu}{i, \, j} d_{j} \equiv \basel{t}{i} \mod \phi(N)$,
for $i = 1$ and $2$, are obtained, then a linear relation of the form
$\sum_{j = 1}^{k} \basel{c}{j} d_{j} \equiv 0 \mod \phi(N)$,
for some integers $\basel{c}{j}$, $1 \leq j \leq k$, not all $0$,
can be found. In addition, if $\rho \, \mid \, \basel{c}{j}$, $1 \leq j \leq k$,
for some integer $\rho \geq 2$, then a relation of the form
$h^{\rho} = 1  \mod N$, for some $h \in \basel{\Integers^{*}}{N}$,
can be found out. Linear relations,
like $\sum_{j = 1}^{k} \basel{c}{j} d_{j} \equiv 0 \mod \phi(N)$,
are called trivial, if it so happens that $\sum_{j = 1}^{k} \basel{c}{j} d_{j} = 0$,
even without applying $\mod \phi(N)$. For quadratic integer sieve,
$\mod 2$ restriction (which can be interpreted as the situation corresponding to $\rho = 2$)
is taken, with a view to improve the efficiency, because
if $g^{2t} = 1 \mod N$, for some integer $t$,
then, with $h = g^{t}$, $(h-1)$ and $(h+1)$ may yield nontrivial
factors of $N$ by $\gcd$.

   The estimation of probability of generating a linear relation in $\basel{d}{j}$, 
   for $1 \leq j \leq k$, does not carry over from elliptic curve method to
   general integer sieve, as the term $(\fieldChar+1)$ plays a pivotal role in
   our estimation of error probabilities of the elliptic curve factorization method.
   As for the primitiveness of the chosen base element $g$, it may be observed that
   the cardinality of $\basel{\Integers^{*}}{N}$ is $\phi(N)$, and among the elements
   of $\basel{\Integers^{*}}{N}$, there are about $\phi\bglb \phi(N) \bgrb$ elements
   that can be primitive (group generator) elements. For multiple base elements,
   the primitiveness constraint may be overcome, but the probability of generating
   a linear relation is less clearly understood.  Subsequently, the merits of
   elliptic curve factorization method are described.

\paragraph{Merits of Elliptic Curve Factorization}
\begin{enumerate}
\item the method is probabilistic polynomial time algorithm under the assumption of uniform probability
of the group orders for random modulus of given size ;  
\item the space requirement is quite small, compared to integer sieve method ;
\item if at least one sample of $k$-smooth group order is realized, then the factorization produces a result ; and
\item it is not necessary to assume that the initial random point for any selected curve is a group generator
\end{enumerate}

 However, diligence must be exercised while exponentiating by a prime number $\basel{\primeInt}{i}$,
 in that the exponentiation may be conducted for at most $\frac{\log(N)}{2\log(\basel{\primeInt}{i})}$
 times, for every positive integer $i \leq k$. The number of curve samples also plays an important role,
 which must be taken in parallel, for each exponentiation by  $\basel{\primeInt}{i}$,  $1 \leq i \leq k$.
 
\section{Conclusion}
The probability analysis for the elliptic curve factorization is presented. 
The method is shown to be a probabilistic polynomial time algorithm, under 
reasonable assumptions on the probability distribution of the group orders
that arise, when restriction to a fixed (but unknown) smaller prime factor
of the modulus integer to be factored is taken. The integer modulus to be
factored is treated as a random variable of fixed size, because it is
an input to the factorization algorithm. The analysis takes into account
the {\em{a priori}} and {\em{a posteriori}} probabilities.
The probability of successful factorization is fairly optimistic.

\end{document}